\documentclass[12pt]{amsart}
\usepackage[margin=1in]{geometry}

\usepackage{url,amsmath,amssymb,mathrsfs}
\usepackage{enumitem}
\usepackage{commath}
\usepackage{physics}
\usepackage{nicefrac}
\usepackage{amsthm}
\usepackage{xfrac}
\usepackage[all]{xy}
\usepackage{graphicx}
\graphicspath{ {./images/} }

\newtheorem{theorem}{Theorem}[section]
\newtheorem{lemma}[theorem]{Lemma}

\newtheorem{proposition}[theorem]{Proposition}

\newtheorem{innercustomgeneric}{\customgenericname}
\providecommand{\customgenericname}{}
\newcommand{\newcustomtheorem}[2]{%
	\newenvironment{#1}[1]
	{%
		\renewcommand\customgenericname{#2}%
		\renewcommand\theinnercustomgeneric{##1}%
		\innercustomgeneric
	}
	{\endinnercustomgeneric}
}
\newcustomtheorem{customthm}{Theorem}
\newcustomtheorem{customlemma}{Lemma}

\theoremstyle{definition}
\newtheorem{example}[theorem]{Example}
\newtheorem{definition}[theorem]{Definition}

\newcommand{\kron}[2]{\left(\frac{#1}{#2}\right)}
\newcommand{\bZ}{\mathbb Z}
\newcommand{\bQ}{\mathbb Q}
\newcommand{\bN}{\mathbb N}
\newcommand{\bR}{\mathbb R}

\newcommand{\Rbar}{\overline{R}}

\newcommand{\Cl}{\textup{Cl}}

\newcommand{\st}{\textup{ s.t. }}

\newcommand{\cO}{\mathcal{O}}
\newcommand{\Irr}{\textup{Irr}}
\newcommand{\term}[1]{\textbf{\textup{#1}}}

\newcommand{\quot}[2]{\large\sfrac{#1}{#2}\normalsize}

\title{Elasticity of Orders with Prime Conductor}
\author{Jared Kettinger \and Grant Moles}
\date{April 2025}

\begin{document}

\maketitle

\begin{abstract}

    Let $R$ be an order in a number field whose conductor ideal $P := (R:\Rbar)$ is prime in the ring of integers $\Rbar$. In this paper, we explore the factorization properties of such orders. Most notably, we give a complete characterization of the elasticity of $R$ in terms of its class group. We conclude with an application to the computation of class groups of certain orders.

\end{abstract}

\section{Introduction}
    Let $R$ be an atomic domain. If every element in $R$ factors uniquely into a product of irreducibles, we call $R$ a \term{unique factorization domain} (UFD). Of course, not every atomic domain exhibits unique factorization. Thus, when determining the structure of the domain $R$, it helps to describe how ``close" $R$ is to having unique factorization. For instance, in 1960, Carlitz described in \cite{carlitz} rings of algebraic integers in which, even though unique factorization may fail, any irreducible factorization of a given element has the same length. Domains exhibiting this factorization property have been investigated at length since; it was Zaks in \cite{zaks1} and \cite{zaks2} who first coined the term \term{half-factorial domain} to describe such a ring.

    Half-factorial domains are certainly interesting, but, in a sense, they only describe the ``best-case scenario" for a domain which fails to exhibit unique factorization. To better understand atomic domains which are not half-factorial and provide a more detailed view of how badly unique factorization fails in these domains, we use the following definition introduced by Valenza \cite{valenza}. Note that in this definition (and throughout this paper), we use the notation $\Irr(R)$ to denote the set of irreducible elements in the ring $R$.
    \begin{definition}
        \label{elasticity}
        Let $R$ be an atomic domain and $\alpha\in R$ a nonzero, nonunit element. The \term{length set} of $\alpha$ in $R$ is
        $$\ell_R(\alpha):=\{k\in\bN|\exists\pi_1,\dots,\pi_k\in\Irr(R)\st \alpha=\pi_1\cdots\pi_k\}.$$
        The \term{elasticity} of $\alpha$ in $R$ is
        $$\rho_R(\alpha):=\sup\left\{\frac{m}{n}\:|\:m,n\in\ell_R(\alpha)\right\}.$$
        The \term{elasticity} of $R$ is
        $$\rho(R):=\sup\{\rho_R(\alpha)|\alpha\in R\backslash(U(R)\cup\{0\})\}.$$
    \end{definition}

    That is, the elasticity of an element $\alpha\in R$ gives the ratio of the longest irreducible factorization of $\alpha$ to its shortest, and the elasticity of $R$ is the supremum of the elasticity of its elements. In particular, an atomic domain $R$ is half-factorial if and only if it has elasticity 1. 
    
    While unique factorization of elements fails for rings of integers in general, we retain unique factorization of (proper) ideals into prime ideals. Thus, for any nonunit $\alpha \in R$, we can use the factorization of the principal ideal $(\alpha) \subseteq R$ to gain information about the factorization of $\alpha$. In particular, we can use the structure of the class group $\Cl(R)$ to inform our knowledge of factorization in $R$. To this end, H. Davenport introduced the following invariant of a finite abelian group in 1966 (\cite{davenport1966midwestern}).

    \begin{definition}
        For a finite abelian group $G$, we call a $G$-sequence $\{g_1,g_2,\dots,g_n\}$ of (not necessarily distinct) elements of $G$ a \term{0-sequence} if $g_1 + g_2 + \dotsm + g_n = 0$. We call a subsequence $\{g_{m_1},g_{m_2},\dots,g_{m_k}\}$ of $\{g_1,g_2,\dots,g_n\}$ which is itself a 0-sequence a \term{0-subsequence}. If $k<n$, that is the 0-subsequence is strictly shorter than the original, we call this a \term{proper 0-subsequence}.
    \end{definition}
    
    \begin{definition}
        For a finite abelian group $G$, we define the \term{Davenport constant} of $G$, denoted $D(G)$, by either of the following equivalent definitions: $D(G) = \text{min}\{n\, | $ any $G$-sequence of length $n$ has a $0$-subsequence$\}$ or $D(G) = \text{max}\{n\, | \, \exists\,$ a $0$-sequence of length $n$ with no proper $0$-subsequence$\}$.
    \end{definition}

    While there has been extensive study of $D(G)$ as a combinatorial object, it was first defined in connection with factorization. In particular, for a ring of integers $R$ and $\alpha \in \text{Irr}(R)$, $D(\Cl(R))$ is a sharp upper bound for the number of prime factors (counting multiplicity) of $(\alpha)$. Using this fact, Valenza in \cite{valenza} was able to determine an upper bound for $\rho(R)$ in terms of $D(\Cl(R))$. Five years later, in 1995, Narkiewicz in \cite{narkiewicz} showed that this bound was actually an equality---giving us the following theorem.

    \begin{theorem}
        Let $\cO_K$ be a ring of integers which is not a UFD, then
        \[
        \rho(\cO_K) = \frac{D(\Cl(\cO_K))}{2}
        \]
    \end{theorem}

    A natural generalization from rings of algebraic integers would be to consider orders in an algebraic number field. Recall that for a number field $K$, an order $R$ in $K$ is a subring (with identity) of the full ring of algebraic integers $\cO_K$ whose quotient field is $K$. Orders in a number field, much like rings of algebraic integers, are Noetherian rings with Krull dimension 1. However, they may fail in general to be integrally closed (and thus may fail to be Dedekind). In fact, if $R$ is an order in $K$, then the integral closure of $R$ is $\Rbar=\cO_K$; then unless the order $R$ is maximal (i.e. is the full ring of integers), $R$ is not integrally closed (and thus is not a UFD).

    It was Halter-Koch in \cite{halter-koch} who made the first major progress in studying factorization in orders in a number field. Specifically, he characterized the orders in a quadratic number field which are half-factorial. More recently, Rago in \cite{rago} gave a full characterization of half-factorial orders in arbitrary number fields. Meanwhile, more general orders whose elasticities match those of their integral closures were studied in \cite{dissertation}. Since 2024, Choi in \cite{choi2024class} and Kettinger in \cite{Kettinger2025} have given results on the elasticity of orders in quadratic rings of integers. Choi's work focused on orders with primary conductor ideal whose integral closures are UFDs, while Kettinger took an element-based approach to orders with prime conductor. Some techniques used in this paper represent a natural, ideal-theoretic analogue of this approach and allow us to greatly expand on the results of \cite{Kettinger2025}.

    In this paper, we will investigate the elasticity of orders in a number field whose conductor ideals are prime in the integral closure (the ring of integers). This will culminate in the main result of this paper, which is presented here.

    \begin{customthm}{\ref{main result}}
        Let $R$ be an order in a number field $K$ with conductor ideal $P$, where $P$ is prime as an ideal of $\Rbar$. Then:
        \begin{enumerate}
            \item If $P$ is principal in $\Rbar$ and the equivalent conditions from Lemma \ref{equiv cond} hold, then
            $$\rho(R)=\frac{D(\Cl(R))+1}{2}.$$
            \item Otherwise,
            $$\rho(R)=\frac{D(\Cl(R))}{2}.$$
        \end{enumerate}
    \end{customthm}

    \section{Preliminaries}

    Before we are able to prove Theorem \ref{main result}, we must first compile the necessary tools. First, we present lemmata which we will reference frequently. The first two come from \cite{conrad}; the third from \cite{picavet}.

    \begin{lemma}
        \label{rel prime to I is invertible}
        Let $R$ be an order in a number field with conductor ideal $I$. If $J$ is an $R$-ideal which is relatively prime to $I$, then $J$ is an invertible ideal of $R$.
    \end{lemma}

    \begin{lemma}
        \label{R ideals factor}
        Let $R$ be an order in a number field with conductor ideal $I$. Any $R$-ideal which is relatively prime to $I$ has unique factorization into prime $R$-ideals relatively prime to $I$. Moreover, all but finitely many prime ideals in $R$ are relatively prime to $I$.
    \end{lemma}

    \begin{lemma}
        \label{prime ideals in classes}
        Let $R$ be an order in a number field. Then every ideal class in $\Cl(R)$ contains infinitely many prime $R$-ideals.
    \end{lemma}

    The next result is a porism to Theorem 3.4 from \cite{radicalconductor} and will give us a lower bound for the elasticity of an order. For convenience, we will present the proof here.

    \begin{proposition}
        \label{elasticity lower bound}
        Let $R$ be an order in a number field. Then $\rho(R)\geq\frac{D(\Cl(R))}{2}$.
    \end{proposition}

    \begin{proof}
        First, note that if $D(\Cl(R))=1$, the result is trivial; then we will assume henceforth that $\Cl(R)$ is nontrivial.

        Let $d=D(\Cl(R))>1$ and let $\{[I_1],[I_2],\dots,[I_d]\}$ be a 0-sequence in $\Cl(R)$ with no proper 0-subsequence. By Lemma \ref{prime ideals in classes}, each of the ideal classes $[I_j]$ contains a prime $R$-ideal $P_j$ which is relatively prime to the conductor ideal $I$ of $R$. Furthermore, each inverse class $[I_j]^{-1}$ similarly contains a prime $R$-ideal $Q_j$ relatively prime to $I$. Now let $\alpha\in R$ be a generator of the principal $R$-ideal $P_1\cdots P_d$; $\beta\in R$ a generator of $Q_1\cdots Q_d$; and for each $1\leq i\leq d$, $\pi_i$ a generator of $P_iQ_i$. Without loss of generality, we will also assume that $\alpha\beta=\pi_1\cdots\pi_d$ (if not, these elements are associates in $R$, and we can absorb the unit into one of the generators).

        Now note that $\alpha$, $\beta$, and each of the $\pi_i$'s must be irreducible in $R$. To see this, suppose that $\alpha=ab$ for some $a,b\in R$; then $\alpha R=(a R)(bR)$. Note that the ideals $\alpha R$, $aR$, and $bR$ are all relatively prime to $I$ by construction, so Lemma \ref{R ideals factor} tells us that these ideals factor uniquely into a product of prime $R$-ideals relatively prime to $I$. Then $\alpha R=P_1\cdots P_d=(aR)(bR)$, meaning that $aR$ factors into some subproduct of the $P_i$'s. However, the sequence $\{[P_1],\dots,[P_d]\}$ has no proper 0-sequence, meaning that either $aR$ is a unit (and factors into an empty product of prime ideals) or that $aR=\alpha R$ (and thus $bR$ is a unit). Then either $a$ or $b$ must be a unit, so $\alpha$ is irreducible in $R$. A similar argument shows that $b$ and each $\pi_i$ must be irreducible. Then $\alpha\beta=\pi_1\cdots\pi_d$, with $\alpha$, $\beta$, and each $\pi_i$ irreducible. Then by definition of elasticity, $\rho(R)\geq \frac{d}{2}=\frac{D(\Cl(R))}{2}$.
    \end{proof}

    From this result (and especially from the main theorem stated earlier), it is clear that the ideal class group of an order will be important when studying elasticity. To give information about the size and structure of this order, we make use of the following result, originally from \cite{narkiewicz}.

    \begin{proposition}
        \label{exact sequence}
        Let $R$ be an order in a number field $K$ with conductor ideal $I$. Then there is an exact sequence
        $$1\to U(R)\to U(\Rbar)\times U(\quot{R}{I})\to U(\quot{\Rbar}{I})\to \Cl(R)\to\Cl(\Rbar)\to 1.$$
        In particular, there is a surjection $\tau:\Cl(R)\to\Cl(\Rbar)$ defined by $\tau([J])=[J\Rbar]$ for every $[J]\in\Cl(R)$. Thus, the class numbers $\abs{\Cl(R)}$ and $|\Cl(\Rbar)|$ are related as follows:
        $$\abs{\Cl(R)}=\abs{\Cl(\Rbar)}\frac{\abs{U(\quot{\Rbar}{I})}}{\abs{U(\quot{R}{I})}\cdot\abs{\quot{U(\Rbar)}{U(R)}}}.$$
    \end{proposition}
    
    As we will see later, the major proofs in this paper will rely heavily on relations between ideals of $R$ and ideals of $\Rbar$. In particular, we will make frequent use of the following lemma that utilizes the surjection $\tau$ from the exact sequence above. Similar results are developed in \cite{choi2024class} and \cite{conrad}; this lemma provides a stronger condition needed for later results.

    \begin{lemma}
        \label{ideals above or below}
        Let $R$ be an order in a number field with conductor ideal $I$. If $J$ is an $R$-ideal relatively prime to $I$, then $J\Rbar\cap R=J$; if $J$ is prime, then $J\Rbar$ is prime as well. If $J$ is an $\Rbar$-ideal relatively prime to $I$, then $(R\cap J)\Rbar=J$; if $J$ is prime, then $R\cap J$ is prime as well.
    \end{lemma}

    \begin{proof}
        First, let $J$ be an $R$-ideal relatively prime to $I$. Since $1\in \Rbar$, $J\subseteq J\Rbar\cap R$ trivially. For the reverse inclusion, let $\alpha=j_1r_1+\dots+j_kr_k$ be an arbitrary element of $J\Rbar\cap R$, with $j_i\in J$ and $r_i\in \Rbar$ for $1\leq i\leq k$. Since $J$ and $I$ are relatively prime as $R$-ideals and $R$ contains $1$ and $\alpha$, there exist $\beta_1,\beta_2\in J$ and $\gamma_1,\gamma_2\in I$ such that $1=\beta_1+\gamma_1$ and $\alpha=\beta_2+\gamma_2$. Then:
        $$\alpha=\alpha(\beta_1+\gamma_1)=(\beta_2+\gamma_2)\beta_1+(j_1r_1+\dots+j_kr_k)\gamma_1=\beta_2\beta_1+\gamma_2\beta_1+j_1r_1\gamma_1+\dots+j_kr_k\gamma_1.$$
        Note that $\beta_2\beta_1\in J^2$, and each remaining term in the above sum is in $JI$. Then $\alpha\in J^2+JI\subseteq J$. Then $J\Rbar\cap R\subseteq J$, and thus $J=J\Rbar\cap R$.

        Now let $J$ be an $\Rbar$-ideal relatively prime to $I$. Trivially, $(R\cap J)\Rbar\subseteq J$. To show the reverse inclusion, let $\alpha+\beta=1$ with $\alpha\in J$ and $\beta\in I$ and note that $\alpha=1-\beta\in R$, so $\alpha\in R\cap J$. Then $1\in R\cap J+I$, so $R\cap J+I=R$ (i.e. $R\cap J$ is relatively prime to $I$ as an $R$-ideal). Then:
        $$J=JR=J(J\cap R+I)=J(J\cap R)+JI\subseteq (J\cap R)\Rbar.$$
        Thus, $J=(J\cap R)\Rbar$.

        If $J$ is a prime $\Rbar$-ideal, then letting $\alpha,\beta\in R$ such that $\alpha\beta\in R\cap J$, note that $\alpha\beta\in J$, so either $\alpha$ or $\beta$ lies in $J$ (and thus in $R\cap J$). Since $R\cap J\neq R$ (in particular, $1\notin R\cap J$), this gives us that $R\cap J$ is a prime $R$-ideal.

        On the other hand, let $J$ be a prime $R$-ideal relatively prime to $I$. Since $J\Rbar\cap R=J\neq R$ (as above), then $J\Rbar\neq \Rbar$. In particular, this means that $J\Rbar\subseteq P$ for some prime $\Rbar$-ideal $P$. Then $P\cap R$ is a prime $R$-ideal containing $J\Rbar\cap R=J$, and since $\dim(R)=1$, $P\cap R=J$. Then $P=(P\cap R)\Rbar=J\Rbar$, so $J\Rbar$ is a prime $\Rbar$-ideal.
    \end{proof}

    When considering the elasticity of an element, we will need to consider the lengths of its longest and shortest irreducible factorizations. The following lemma will help in determining the length of the longest factorization; the theorem after will help in determining the length of the shortest factorization.

    \begin{lemma}
        \label{no prime factors}
        Let $R$ be an order in a number field $K$ with conductor ideal $I$. Let $\alpha\in R$ with $\rho_R(\alpha)>1$ and suppose that $\alpha=\beta\pi$ for some $\beta\in \Rbar$ and $\pi\in R$, with $\pi$ a prime element of $\Rbar$ relatively prime to $I$. Then $\beta\in R$ and $\rho_R(\beta)\geq \rho_R(\alpha)$.
    \end{lemma}

    \begin{proof}
        Since $\pi\in R$ is relatively prime to $I$, then reducing modulo $I$ gives $\beta+I=(\alpha+I)(\pi+I)^{-1}\in \quot{R}{I}$ (since $\pi+I\in \quot{R}{I}\cap U(\quot{\Rbar}{I})$, with $\quot{\Rbar}{I}$ integral over $\quot{R}{I}$). Then $\beta\in R$. If $\beta\in U(R)$, then $\alpha$ is an associate of $\pi$, an irreducible element of $R$. Then $\alpha$ is also irreducible in $R$ and thus $\rho_R(\alpha)=1$, a contradiction. Then $\beta\notin U(R)$ (and thus we can consider its elasticity).

        Now say $\alpha=\pi_1\cdots\pi_m=\tau_1\cdots\tau_n$ with $\pi_i,\tau_j\in \Irr(R)$ for $1\leq i\leq m$, $1\leq j\leq n$ and $n>m$. Since $\pi|\alpha$ and $\pi$ is prime in $\Rbar$, necessarily we have some $\pi_i$ and $\tau_j$ which are associates of $\pi$. Without loss of generality assume that $\pi=\pi_1$ and $\tau_1=u\pi$ for some $u\in U(\Rbar)$. Again, since $\pi$ is relatively prime to $I$, we get that $u\in U(R)$. Then $\beta=\frac{\alpha}{\pi}=\pi_2\cdots\pi_m=(u\tau_2)\cdots \tau_n$. Then for any pair of factorizations of lengths $m$ and $n$ of $\alpha$, $\beta$ has a pair of factorizations of lengths $m-1$ and $n-1$, with $\frac{n}{m}<\frac{n-1}{m-1}$. Therefore, $\rho_R(\alpha)\leq \rho_R(\beta)$.
    \end{proof}

    Note the utility of this lemma. When trying to find the elasticity of the order $R$, we will be trying to find elements of large elasticity. The lemma ensures that if $\alpha\in R$ has any irreducible factors which are prime in $\Rbar$ and relatively prime to $I$, then we can remove those prime factors to produce an element with larger elasticity. When considering an element $\alpha\in R$ of large elasticity, we can therefore assume that $\alpha\Rbar$ has no principal prime ideal factors in $\Rbar$ except for possibly those dividing $I$.

    \begin{theorem}\label{IrrBound}
        \label{number of prime ideal factors}
        Let $R$ be an order in a number field $K$ with prime conductor ideal $P$. Then if $\alpha\in \Irr(R)$, the principal ideal $\alpha\Rbar$ factors into at most $D(\Cl(R))$ prime ideals in $\Rbar$.
    \end{theorem}

    \begin{proof}
        Let $\alpha\in\Irr(R)$. First, assume that $\alpha$ is relatively prime to $P$. By Lemma \ref{R ideals factor}, this means that the ideal $\alpha R$ uniquely factors into a product of prime $R$-ideals $\alpha R=Q_1\cdots Q_k$. Since each $Q_i$ is relatively prime to $P$, these prime ideals are invertible by Lemma \ref{rel prime to I is invertible} and thus represent ideal classes in $\Cl(R)$. Since $\alpha$ is irreducible in $R$, the sequence $\{[Q_1],\dots,[Q_k]\}$ must be a 0-sequence in $\Cl(R)$ with no 0-subsequence, and thus $k\leq D(\Cl(R))$. Then since $\alpha \Rbar=(Q_1\Rbar)\cdots(Q_k\Rbar)$ (with each $Q_i\Rbar$ a prime $\Rbar$-ideal by Lemma \ref{ideals above or below}), $\alpha\Rbar$ factors into at most $D(\Cl(R))$ prime ideals in $\Rbar$.
        
        Now drop the assumption that $\alpha$ is relatively prime to $P$ and write $\alpha\Rbar=P^a P_1\cdots P_k$, with each $P_i$ a prime ideal of $\Rbar$ relatively prime to $P$ and $a\in\bN$. Assume first that $a=1$, and suppose that some subproduct of $P_1\cdots P_k$ is a principal ideal generated by an element of $R$. That is, after possible rearranging, there is some $1\leq j\leq k$ such that $P_1\cdots P_j=\beta \Rbar$ with $\beta\in R$. Then necessarily $PP_{j+1}\cdots P_k$ is principal as well and generated by some $\gamma\in R$ (since $\gamma\in P$, the conductor ideal of $R$). Then for some unit $u\in U(\Rbar)$, $\alpha=\beta(u\gamma)$, with $\beta\in R$ and $u\gamma\in P\subseteq R$. Since neither $\beta$ nor $u\gamma$ is a unit (as a generator of a proper ideal of $\Rbar$), this means that $\alpha$ is not irreducible in $R$, contrary to its definition. Then no subproduct of $P_1\cdots P_k$ can be principal and generated by an element of $R$. Now consider the sequence $\{[R\cap P_1],\dots,[R\cap P_k]\}$ in $\Cl(R)$, and note by Lemma \ref{ideals above or below} that this cannot have a 0-subsequence (otherwise, $(R\cap P_1)\cdots(R\cap P_j)=\beta R$ would imply that $P_1\cdots P_j=(R\cap P_1)\Rbar\cdots(R\cap P_j)\Rbar=\beta\Rbar$). Then $k<D(\Cl(R))$, so $\alpha\Rbar$ factors into at most $D(\Cl(R))$ prime ideals in $\Rbar$.

        Finally, let $\alpha\Rbar=P^aP_1\cdots P_k$ with $a\geq 2$, and let $Q$ be a prime ideal of $\Rbar$ relatively prime to $P$ such that $[Q]=[P]$. Then consider the sequence $\{[R\cap Q]^{a-1},[R\cap P_1],\dots,[R\cap P_k]\}$ of length $a+k-1$ in $\Cl(R)$ (here, $[R\cap Q]^{a-1}$ refers to $a-1$ copies of the element $[R\cap Q]$). Suppose that this sequence had a 0-subsequence. If that 0-subsequence were of the form $\{[R\cap P_1],\dots,[R\cap P_j]\}$ for some $1\leq j\leq k$, then $\alpha$ would factor as $\alpha=\beta\gamma$, with $\beta\in R$ a generator of $P_1\cdots P_j$ and $\gamma$ a generator of $P^aP_{j+1}\cdots P_k\in P\subseteq R$. Much as in the $a=1$ case, this would contradict the irreducibility of $\alpha$ in $R$. Otherwise, the 0-sequence must be of the form $\{[R\cap Q]^b,[R\cap P_1],\dots,[R\cap P_j]\}$ for some $1\leq b<a$ and $0\leq j\leq k$. In this case, $\{[Q]^b,[P_1],\dots,[P_j]\}$ must be a 0-sequence in $\Cl(\Rbar)$ (as the image of a 0-sequence in $\Cl(R)$ under the homomorphism $\tau$ from Proposition \ref{exact sequence}). Since $[P]=[Q]$, this means that $P^bP_1\cdots P_j$ is a principal ideal. Then $P^{a-b}P_{j+1}\cdots P_k$ is also principal; let $\beta$ and $\gamma$, respectively, be the generators of these ideals. Then for some $u\in U(\Rbar)$, $\alpha=\beta(u\gamma)$. As both $\beta$ and $\gamma$ are elements of $P\subseteq R$, this produces a non-trivial factorization of $\alpha$ in $R$, a contradiction. Thus, the sequence $\{[Q]^{a-1},[P_1],\dots,[P_k]\}$ must in fact have no 0-subsequence, so $a+k-1<D(\Cl(R))$. Then $a+k\leq D(\Cl(R))$, so $\alpha\Rbar$ again factors into at most $D(\Cl(R))$ prime ideals in $\Rbar$.
    \end{proof}

    The proof of Proposition \ref{elasticity lower bound} shows that the upper bound given in Theorem \ref{number of prime ideal factors} is in fact sharp for any order $R$. Note that in this proof, the ideal generated by the irreducible elements $\alpha,\beta\in R$ factored into exactly $D(\Cl(R))$ prime $\Rbar$-ideals. Moreover, such a factorization is possible regardless of the conductor ideal of the order.




    \section{Main Result}

    Before presenting the proof of Theorem \ref{main result}, we introduce a final lemma to simplify the theorem statement and to provide an insight into the conditions which affect the elasticity of the order.

    \begin{lemma}
        \label{equiv cond}
        Let $R$ be an order in a number field whose conductor ideal $P$ is a principal prime ideal of $\Rbar$; that is, there exists a prime element $\pi\in R$ such that $P=\pi\Rbar$. The following are equivalent:
        \begin{enumerate}
            \item There exists $\beta\in\Rbar$ such that no nonunit in $R$ divides $\beta$ and $\beta\Rbar$ factors into exactly $D(\Cl(R))-1$ prime ideals of $\Rbar$.
            \item There exists an ideal $A$ of $R$ such that $A\Rbar$ is a principal $\Rbar$-ideal and $[A]$ is contained in a 0-sequence in $\Cl(R)$ of length $D(\Cl(R))$ with no proper 0-subsequence.
        \end{enumerate}
    \end{lemma}

    \begin{proof}
        First, assume that there exists an element $\beta\in\Rbar$ as in Condition (1). Note that $\beta\notin P$; otherwise, $\pi|\beta$ with $\pi\in R$, contradicting the definition of $\beta$. Then letting $d=D(\Cl(R))$, we can factor $\beta\Rbar=P_1\cdots P_{d-1}$ with each $P_{i}$ a prime ideal of $\Rbar$ relatively prime to $P$. Then consider the sequence $\{[R\cap P_1],\dots,[R\cap P_{d-1}]\}$ in $\Cl(R)$. If this sequence had a 0-subsequence, say $\{[R\cap P_1],\dots,[R\cap P_j]\}$, then there would exist $\gamma\in R$ such that $\gamma R=(R\cap P_1)\cdots(R\cap P_j)$ and thus $\gamma\Rbar=(R\cap P_1)\Rbar\cdots(R\cap P_j)\Rbar=P_1\cdots P_j$ (by Lemma \ref{ideals above or below}). Since such $\gamma\in R$ would be a nonunit divisor of $\beta$, no such $\gamma$ can exist, and thus the sequence $\{[R\cap P_1],\dots, [R\cap P_k]\}$ in $\Cl(R)$ can have no 0-subsequence. 
        
        Now letting $[A]=[\beta R]^{-1}\in \Cl(R)$ with $A$ chosen to be an ideal of $R$ relatively prime to $P$ (such $A$ must exist by Lemmas \ref{R ideals factor} and \ref{prime ideals in classes}), we have that $[A]$ is contained in the 0-sequence $\{[R\cap P_1],\dots,[R\cap P_{d-1}],[A]\}$ of length $d$. Note that this 0-sequence has no proper 0-subsequence since $\{[R\cap P_1],\dots,[R\cap P_{d-1}]$ has no 0-subsequence, and $A\Rbar$ must be principal since $\beta\Rbar$ and $(\beta\Rbar)(A\Rbar)$ are both principal. Then $(1)\implies (2)$.

        Now assume that there exists an ideal $A$ of $R$ as in Condition (2), and let $\{[A],[P_1],\dots,[P_{d-1}]\}$ be a 0-sequence of length $d=D(\Cl(R))$ with no proper 0-sequence. As usual, we will assume that the $P_i$'s are prime ideals of $R$ relatively prime to $P$. Letting $\alpha R=AP_1\cdots P_{d-1}$, note that $\alpha\Rbar=A\Rbar(P_1\Rbar)\cdots (P_{d-1}\Rbar)$, with $\alpha\Rbar$ and $A\Rbar$ both principal. Then for some $\beta\in\Rbar$, $\beta\Rbar=(P_1\Rbar)\cdots (P_{d-1}\Rbar)$ (that is, $\beta$ factors into exactly $D(\Cl(R))-1$ prime ideals of $\Rbar$). Furthermore, if there exists a nonunit $\gamma\in R$ such that $\gamma|\beta$, then since $\beta$ (and thus $\gamma$) is relatively prime to $P$, we can use Lemma \ref{R ideals factor} to write $\gamma R=Q_1\cdots Q_k$, with each $Q_i$ a prime ideal of $R$ relatively prime to $P$. Then $\gamma \Rbar=(Q_1\Rbar)\cdots (Q_k\Rbar)$, and since $\gamma|\beta$, necessarily (after possibly rearranging the $P_i$'s) $1\leq k\leq d-1$ and $Q_i=P_i$ for each $1\leq i\leq k$. Then $\{[Q_1],\dots,[Q_k]\}$ is a proper 0-subsequence of $\{[A],[P_1],\dots,[P_{d-1}]\}$, a contradiction. Thus, $\beta$ can have no nonunit divisors lying in $R$. Therefore, $(2)\implies (1)$.
    \end{proof}

    We now restate and prove the main result of this paper.

    \begin{theorem}
    \label{main result}
        Let $R$ be an order in a number field $K$ with conductor ideal $P$, where $P$ is prime as an ideal of $\Rbar$. Then:
        \begin{enumerate}
            \item If $P$ is principal in $\Rbar$ and either of the equivalent conditions from Lemma \ref{equiv cond} hold, then
            $$\rho(R)=\frac{D(\Cl(R))+1}{2}.$$
            \item Otherwise,
            $$\rho(R)=\frac{D(\Cl(R))}{2}.$$
        \end{enumerate}
        In particular, if $\Cl(\Rbar)$ is trivial, then case 1 holds (if $\Cl(R)=1$ as well, then $\rho(R)=1$). If $\Cl(R)\cong \Cl(\Rbar)\not\cong\{1\}$, then case 2 holds.
    \end{theorem}

    \begin{proof}
        For ease of notation, we will denote $d=D(\Cl(R))$ throughout this proof. First, note that if $\abs{\Cl(R)}=1$, then $|\Cl(\Rbar)|=1$ as well. As an immediate consequence of Proposition 2.4.5, Corollary 2.4.7, and Theorem 3.1.7 from \cite{dissertation}, $\rho(R)=\rho(\Rbar)=1$. Note that this agrees with case 1 above, since any ideal $A$ of $R$ relatively prime to $P$ automatically satisfies Condition 2 from Lemma \ref{equiv cond}.

        When determining the elasticity of $R$, we need to consider elements of $R$ with large elasticity. By Lemma \ref{no prime factors} and the discussion that followed, we need only consider elements $\alpha\in R$ such that $\alpha\Rbar$ has no principal prime ideal factors in $\Rbar$ other than possibly $P$ itself. Then letting $\alpha\in R$ be a nonzero, nonunit, we will henceforth assume that $\alpha\Rbar$ has no principal prime ideal factors in $\Rbar$ except possibly $P$.

        Suppose that $P$ is non-principal in $\Rbar$. Since $P$ is non-principal, we want to show that $\rho(R)=\frac{d}{2}$. By Proposition \ref{elasticity lower bound}, we have that $\rho(R)\geq \frac{d}{2}$; we must show the reverse inequality. Let $\alpha\in R$ and write $\alpha=\pi_1\cdots \pi_m=\tau_1\cdots\tau_n$ with $\pi_i,\tau_j\in\Irr(R)$ and $n\geq m$. By the above discussion, we may assume that none of the $\pi_i$'s or $\tau_j$'s are prime elements of $\Rbar$ (since $P$ is non-principal). Then each principal ideal $\pi_i\Rbar$ and $\tau_j\Rbar$ must factor into at least $2$ prime $\Rbar$-ideals (since each $\pi_i$ and $\tau_j$ is a nonunit in $R$ and must thus be a nonunit in $\Rbar$ by integrality). By Theorem \ref{number of prime ideal factors}, each $\pi_i\Rbar$ and $\tau_j \Rbar$ must factor into at most $d$ prime $\Rbar$-ideals. Then if $\alpha\Rbar=P_1\cdots P_k$ is the prime ideal factorization of $\alpha\Rbar$ in $\Rbar$, we have $2n\leq k\leq dm$ and thus $\frac{n}{m}\leq \frac{d}{2}$. Then the elasticity $\rho_R(\alpha)$ of any nonzero, nonunit element $\alpha\in R$ must be at most $\frac{d}{2}$, so $\rho(R)\leq \frac{d}{2}$. Then since we have shown both inequalities, $\rho(R)=\frac{d}{2}$.

        Now suppose that $P$ is a principal ideal of $\Rbar$ generated by $\pi\in\Rbar$; since $P$ is the conductor ideal of $R$, necessarily $\pi\in R$. Now when considering irreducible factorizations of a nonzero, nonunit $\alpha\in R$, we can no longer assume that $\alpha\Rbar$ has no principal prime ideal factors, since $P$ may divide $\alpha\Rbar$. Then given an irreducible factor $\tau_i$ of $\alpha$, we can no longer guarantee that $\tau_i\Rbar$ factors into at least 2 prime ideal factors. However, note that since $\pi$ conducts every element of $\Rbar$ into $R$, then $\pi^2$ cannot divide any irreducible element of $R$ (if $\beta=\pi^2\gamma$ with $\gamma\in\Rbar$, then $\beta=\pi(\pi\gamma)$ with both $\pi$ and $\pi\gamma$ nonunits in $R$). Then again writing $\alpha=\pi_1\cdots\pi_m=\tau_1\cdots \tau_n$ with each $\pi_i,\tau_j\in\Irr(R)$ and $n\geq m$, if $\pi^r$ is the largest power of $\pi$ dividing $\alpha$, then necessarily $r\leq m\leq n$. Then at most $r$ of the $\tau_j$'s will be prime elements of $\Rbar$; the rest must have at least two prime ideal factors. Therefore, if $\alpha\Rbar=P_1\cdots P_k$ is the prime factorization of $\alpha\Rbar$ in $\Rbar$, then $k\geq r+2(n-r)=2n-r$.

        Let $\beta\in\Rbar$ such that no nonunit in $R$ divides $\beta$. Since $\beta$ itself cannot lie in $R$ (unless $\beta\in U(R)$), then $\beta\Rbar$ must factor into at most $d-1$ prime ideal factors in $\Rbar$ (otherwise, as seen previously, $\beta$ would have a nonunit factor in $R$). Moreover, this upper bound $d-1$ is achievable for some $\beta\in\Rbar$ if and only if we are in case 1 above. In other words, if we define $a\in\bN$ to be minimal such that for some $\beta\in\Rbar$ with no nonunit divisors in $R$, $\beta\Rbar$ factors into exactly $d-a$ primes in $\Rbar$, then $a=1$ if and only if we are in case 1.

        Now considering the irreducible factorization $\alpha=\pi_1\cdots \pi_m$, suppose that for some $1\leq i\leq m$, $\pi_i\in P$. Then necessarily $\pi_i=\pi\beta$ for some $\beta\in \Rbar$. Moreover, if $\beta=\gamma_1\gamma_2$ with $\gamma_1$ a nonunit element of $R$, then $\pi_i=(\pi\gamma_2)\gamma_1$, contradicting the irreducibility of $\pi_i$ in $R$. Then $\beta$ must have no nonunit factors in $R$, and thus by the discussion above, $\pi_i\Rbar$ can factor into at most $d-a+1$ prime ideals in $\Rbar$. Since each $\pi_i$ lying in $P$ can have at most $d-a+1$ prime ideal factors and each $\pi_i$ relatively prime to $P$ can have at most $d$ prime ideal factors, then $r(d-a+1)+d(m-r)\geq k\geq 2n-r$. Rearranging:
        $$\frac{n}{m}\leq \frac{d-\frac{r}{m}(a-2)}{2}.$$

        If $a\geq 2$ (that is, if we are in case 2), then $a-2\geq 0$. Then the fraction on the right-hand side of the above inequality is maximized when $r=0$ (that is, when $\alpha\notin P$). Then in case 2, the elasticity of any element of $R$ is bounded above by $\frac{d}{2}$, so $\rho(R)\leq \frac{d}{2}$. With Proposition \ref{elasticity lower bound}, this gives us that $\rho(R)=\frac{d}{2}$ in case 2.

        If $a=1$ (that is, if we are in case 1), then $a-2=-1$. Then the fraction on the right-hand side of the above inequality is maximized when $r=m$ (that is, when each $\pi_i\in P$). Thus, the elasticity of any element of $R$ is bounded above by $\frac{d+1}{2}$, so $\rho(R)\leq \frac{d+1}{2}$ in case 1. All that remains to show is that in this case, there exists an element $\alpha\in R$ with elasticity exactly $\frac{d+1}{2}$.

        Suppose that we are in case 1 and $\abs{\Cl(R)}\geq 2$ (the case when $\abs{\Cl(R)}=1$ was handled earlier). Then using Condition 2 from Lemma \ref{equiv cond}, there exists an ideal $A$ of $R$ such that $A\Rbar$ is a principal $\Rbar$ ideal and there exists a 0-sequence $\{[A],[A_1],\dots,[A_{d-1}]\}$ in $\Cl(R)$ with no proper 0-subsequence. For each $1\leq i\leq d-1$, let $P_i\in[A_i]$ and $Q_i\in [A_i]^{-1}$ be prime ideals relatively prime to $P$. Since $A\Rbar$ and $(A\Rbar)(P_1\Rbar)\cdots(P_{d-1}\Rbar)$ are both principal in $\Rbar$, then there is some $\alpha\in\Rbar$ such that $\alpha\Rbar=(P_1\Rbar)\cdots(P_{d-1}\Rbar)$. Similarly, there exist $\beta\in \Rbar$ and $\pi_i\in\Rbar$ for $1\leq i\leq d-1$ such that $\beta\Rbar=(Q_1\Rbar)\cdots(Q_{d-1}\Rbar)$ and $\pi_i=(P_i\Rbar)(Q_i\Rbar)$. Assuming without loss of generality that $\alpha\beta=\pi_1\cdots\pi_{d-1}$, we have that $(\pi\alpha)(\pi\beta)=\pi^2\pi_1\cdots\pi_{d-1}$. Since $\pi$ is irreducible in $\Rbar$, it is also irreducible in $R$. Since $D(\Cl(R))>1$, then each $P_i$ and $Q_i$ must be a non-principal ideal of $R$, so each $\pi_i$ must be an irreducible element of $R$. Finally, note that if $\pi\alpha=rs$ with $r$ and $s$ nonunits in $R$, then without loss of generality, $\pi|r$ in $\Rbar$. Then $\alpha=\frac{r}{\pi}\cdot s$, and thus $\alpha$ has a nonunit divisor $s$ which lies in $R$, contrary to assumption. Then $\pi\alpha$ (and by a similar argument, $\pi\beta$) is irreducible in $R$ as well. Then we have constructed an element in $R$ with elasticity at least $\frac{d+1}{2}$; since this is also an upper bound for the elasticity of any element of $R$, this element must have elasticity exactly $\frac{d+1}{2}$. Then in case 1, $\rho(R)=\frac{d+1}{2}$, completing the proof.
    \end{proof}

    Despite the fact that an order fails to be Dedekind in general, Theorem \ref{main result} shows us that---in the case where the conductor ideal is prime in the integral closure---its elasticity remains completely determined by the structure of its class group. The theorem also enjoys a nice symmetry with Narkiewicz's result on elasticity of rings of integers. Recall that for a ring of integers $R$, $\rho(R) = 1 = \frac{D(\Cl(R))+1}{2}$ if $R$ is a UFD (the analogue to case $1$), and $\rho(R) = \frac{D(\Cl(R))}{2}$ otherwise.

\begin{example}
    Consider the quadratic ring of integers $\Rbar = \bZ[\frac{1+\sqrt{-7}}{2}]:=\bZ[\omega]$ and the order $R = \bZ[5 \cdot \frac{1+\sqrt{-7}}{2}] = \bZ[5\cdot \omega]$ with prime conductor $5\Rbar$ and $\Cl(R) \cong \bZ/6\bZ$.  As $\Cl(\Rbar)$ is trivial, we are in case $1$. By Theorem \ref{main result}, we must have $\rho(\bZ[5\cdot \omega]) = \frac{D(\Cl(R))+1}{2} = \frac{6+1}{2} = \frac{7}{2}$. In particular, $800 = (5\omega^5)(5 \bar{\omega}^5) = 5^2(\omega \bar{\omega})^5 =  5^2\cdot 2^5$ achieves elasticity $\frac{7}{2}$.
\end{example}

A natural question to ask at this point is whether it is possible for an order to land in case 1 without the ring of integers being a UFD. The following example addresses this question.

\begin{example}
    Let $R=\bZ[17\sqrt{10}]$, the index 17 order in $\bR=\bZ[\sqrt{10}]$. The rational prime 17 is inert in $\bR$, so the conductor ideal $I=17\Rbar$ is prime in $\Rbar$. $\bR$ is an HFD (i.e. $|\Cl(\Rbar)|=2$), and the fundamental unit in $U(\Rbar)$ is $u=3-\sqrt{10}$, and the minimal power of this unit lying in $R$ is $u^9$. Then using Proposition \ref{exact sequence}, we have
    $$\abs{\Cl(R)}=\abs{\Cl(\Rbar)}\frac{\abs{U(\quot{\overline{R}}{I})}}{\abs{U(\quot{R}{I})}\cdot\abs{\quot{U(\Rbar)}{U(R)}}}=2\cdot\frac{288}{16\cdot 9}=4.$$
    Then $\Cl(R)$ must be isomorphic to either $\bZ_4$ or $\bZ_2\times\bZ_2$. In $\Rbar$, the rational prime 2 factors as $2\Rbar=(2,\sqrt{10})^2$; as there is no element in $\Rbar$ of norm $\pm 2$, the prime ideal $(2,\sqrt{10})=2\bZ+\sqrt{10}\bZ$ is non-principal. Since $R=\bZ+17\sqrt{10}\bZ$, then $R\cap P=2\bZ+17\sqrt{10}\bZ=(2,17\sqrt{10})$. Finally, note that $(2,17\sqrt{10})^2=(4,34\sqrt{10},2890)=2R$, a principal ideal. In other words, $R\cap P$ is an invertible ideal (as it is relatively prime to $I$) of $R$ which does not lie in $\ker(\tau)$ (as $(R\cap P)\Rbar=P$ is not principal in $\Rbar$), but whose order in $\Cl(R)$ is 2. If $\Cl(R)\cong\bZ_4$, then any element not in $\ker(\tau)\cong\bZ_2$ would necessarily have order 4. Then $\Cl(R)$ must be isomorphic to $\bZ_2\times\bZ_2$. Now note that $\{(1,0),(0,1),(1,1)\}$ is a 0-sequence of length $D(\Cl(R))=3$ in $\bZ_2\times\bZ_2$, and that the non-identity element in $\ker(\tau)\cong\bZ_2$ must be in this 0-sequence (as every non-identity element of $\Cl(R)$ is in this sequence). Then $R$ satisfies the condition for case 1, and thus $\rho(R)=\frac{D(\Cl(R))+1}{2}=2$.
\end{example}

The following example demonstrates what can go wrong if the conductor ideal is not prime.

\begin{example}
    Consider the ring of integers $\Rbar = \bZ[i]$ and the order $R = \bZ[5i]$ with conductor $I = 5\bZ[i]$. In $\Rbar$, $I$ factors into prime ideals as $I = (1+2i)(1-2i)$. This gives rise to the irreducible factorizations $(5(1+2i)^n)(5(1-2i)^n) = 5^2((1+2i)(1-2i))^n = 5^{n+2}$ in $R$. Hence, $\rho(R) \geq \frac{n+2}{2}$ for all $n \in \bN$, so we have $\rho(R) = \infty$. 
\end{example}

For a characterization of orders with finite elasticity, see \cite{halter1995elasticity}. The following example tells us that there exists a wealth of orders satisfying case $2$.

\begin{example}
    Let $\Rbar$ be a ring of integers with non-trivial class group, and let $P$ be a non-principal prime ideal of $\Rbar$ with inertial degree $f(P:p) \geq 2$ (i.e. $\abs{\sfrac{\Rbar}{P}}$ is not a rational prime). If we define $R = \bZ + P$, then by \cite{conrad}, $(R:\Rbar) = P$. Therefore, by Theorem \ref{main result}, we have $\rho(R) = \frac{D(\Cl(R))}{2}$. 
\end{example}

As a specific example of an order of this form which has been of prior interest, let $\alpha$ be a root of $x^3+4x-1$, $P=(3,2+2\alpha+\alpha^2)$ in $\Rbar=\bZ[\alpha]$, and $R=\bZ+P$. Then $\rho(R)=\frac{D(\Cl(R))}{2}=1$. The details for this order can be found in \cite{radicalconductor}. The following example gives a more interesting order of the form described here. In particular, this order is one whose integral closure is not an HFD and whose elasticity cannot be found using the results found in \cite{radicalconductor}.

\begin{example}
    Let $\Rbar=\bZ[\alpha]$, the ring of integers in $K=\bQ[\alpha]$ with $\alpha$ a root of $x^4-45x^2+1$. Then $\Cl(\Rbar)\cong \bZ_3$ and $P=(11,1+5\alpha+\alpha^2)$ is a non-principal prime ideal of $\Rbar$ with inertial degree 2. Letting $R=\bZ+P$, we find that $\abs{U(\sfrac{\Rbar}{P})}=120$ and $\abs{U(\sfrac{R}{P})}=10$. Utilizing SageMath, we find that $U(\Rbar)$ has three fundamental units: $\alpha$, $-48-322\alpha+7\alpha^3$, and $409-2743\alpha-9\alpha^2+61$. Using these units, we again use SageMath to determine that $\quot{U(\Rbar)}{U(R)}\cong \bZ_6$ (generated by $\alpha U(R)$). Thus,
    $$\abs{\Cl(R)}=\abs{\Cl(\Rbar)}\frac{\abs{U(\quot{\Rbar}{I})}}{\abs{U(\quot{R}{I})}\cdot\abs{\quot{U(\Rbar)}{U(R)}}}=3\cdot\frac{120}{10\cdot 6}=6.$$
    As there is only one abelian group of order 6, we have $\Cl(R)\cong \bZ_6$, so $\rho(R)=\frac{D(\Cl(R))}{2}=3$.
\end{example}

\section{The Structure of $\Cl(R)$}

Recall that when determining the structure of the class group of a ring of algebraic integers, we rely on powerful tools, such as the Minkowski bound, to do so. Determining the structure of the class group of a more general order in a number field is not as well understood. That said, tools such as Proposition \ref{exact sequence} can be used to our advantage.

This result is certainly useful, but it falls short of telling us the structure of the group $\Cl(R)$; rather, it gives us homomorphisms (in particular, the surjection $\tau$) which describe $\Cl(R)$ and a way of determining the order of the group. Since the major results of this paper depend on knowing the structure of the class group (not just its order), it may be useful to have additional methods of determining this structure. In this section, we utilize the following result together with Theorem \ref{main result} to develop methods for determining the structure of $\Cl(R)$.

\begin{theorem}
    \label{cylic class group}
    Let $R$ be an order in a number field with prime conductor ideal $P$ such that $|\Cl(\Rbar)|=1$ (i.e. $\Rbar$ is a UFD). Then $\Cl(R)$ is a cyclic group of order $$\abs{\Cl(R)}=\frac{\abs{U(\quot{\Rbar}{P})}}{\abs{U(\quot{R}{P})}\cdot\abs{\quot{U(\Rbar)}{U(R)}}}.$$
\end{theorem}


\begin{proof}
    The order of $\Cl(R)$ presented here follows immediately from Proposition \ref{exact sequence}. All that remains to show is that $\Cl(R)$ is indeed a cyclic group. Letting $\tau:\Cl(R)\to\Cl(\Rbar)$ be the surjective map from the exact sequence in Proposition \ref{exact sequence}, it follows from the First Isomorphism Theorem that $$\quot{\Cl(R)}{\ker(\tau)}\cong\Cl(\Rbar)=\{1\},$$ and thus $\Cl(R)\cong \ker(\tau)$. By taking advantage of the exact sequence in Proposition \ref{exact sequence}, it also follows that $$\ker(\tau)\cong \quot{G}{H},$$ where $$G=\quot{U(\sfrac{\Rbar}{P})}{U(\sfrac{R}{P})}$$ and $H$ is the subgroup of $G$ consisting of the cosets modulo $P$ which contain units in $U(\Rbar)$ (for more details of this proof, see the proofs of Theorems 3.1.1 and 3.2.1 from \cite{thesis}). Since $P$ is a prime ideal, then $\quot{\Rbar}{P}$ is a finite field; then $U(\quot{\Rbar}{P})$ is a cyclic group. Then $G$ (as a quotient of a cyclic group) is also cyclic, and thus $\Cl(R)\cong \ker(\tau)$ is cyclic (as a quotient of a cyclic group).
\end{proof}

\begin{example}
    Let $\alpha=\frac{1+\sqrt{641}}{2}$ and $R=\bZ[449\alpha]$, the index 449 order in the UFD $\Rbar=\bZ[\alpha]$. Since 449 is an inert prime in $\Rbar$, then $R$ has prime conductor ideal $P=449\Rbar$. Then $\abs{U(\sfrac{\Rbar}{P})}=449^2-1=201600$ and $\abs{U(\sfrac{R}{I})}=449-1=448$. Letting $u$ be the fundamental unit in $\Rbar$, we find that the minimal power of $u$ lying in $R$ is $u^{3}$. Then
    $$\abs{\Cl(R)}=\frac{201600}{448\cdot 3}=150.$$
    There are two abelian groups of order 150: $\bZ_{150}$ and $\bZ_{5}\times\bZ_{30}$. The theorem guarantees that $\Cl(R)\cong\bZ_{150}$, the cyclic group of order 150. Moreover, since $\Rbar$ is a UFD, we can conclude from Theorem \ref{main result} that $\rho(R)=\frac{D(\Cl(R))+1}{2}=\frac{151}{2}$. The details for this example were drawn from \cite{lmfdb}, the table at \cite{quadla}, and SageMath.
\end{example}


The following result gives a concrete example of how one might determine the structure of $\Cl(R)$ using Proposition \ref{exact sequence}, Theorem \ref{main result}, and specific known factorizations in $R$. 

\begin{example}
    Let $R=\bZ[11\cdot \sqrt{-14}]$, the index 11 order in the quadratic number field $K=\bQ[\sqrt{-14}]$. Since $\kron{-14}{11}=-1$, the conductor ideal $I=11\Rbar$ of $R$ is a prime ideal. The class group of $\Rbar=\bZ[\sqrt{-14}]$ is $\Cl(\Rbar)\cong \bZ_4$. Since $K$ is a non-real quadratic number field with $d=-14<-3$, then $U(\Rbar)=U(R)=\{\pm 1\}$. Then letting $\tau:\Cl(R)\to \Cl(\Rbar)$ be the surjection defined in Proposition \ref{exact sequence}, we have:
    $$\abs{\ker(\tau)}=\frac{\abs{U(\quot{\Rbar}{I})}}{\abs{U(\quot{R}{I})}}=\frac{11^2-1}{11-1}=12.$$
    Since $I$ is a prime ideal, $\quot{\Rbar}{I}$ is a finite field and thus has cyclic unit group; then $\ker(\tau)\cong\bZ_{12}$. Then we have:
    $$\quot{\Cl(R)}{\ker(\tau)}\cong \Cl(\Rbar)\cong \bZ_4,$$
    which leaves only three possibilities for $\Cl(R)$: it must be isomorphic to either $\bZ_4\times\bZ_{12}$, $\bZ_2\times\bZ_{24}$, or $\bZ_{48}$. Note that this also leaves only three possibilities for $D(\Cl(R))$: 15, 25, or 48. Now letting $P=(19,9+\sqrt{-14})$, we note that $[P]$ is a generator of the class group $\Cl(\Rbar)$, so $P^4=(325+42\sqrt{-14})$ is irreducible in $\Rbar$. Moreover, $(325+42\sqrt{-14})^{11}$ is an element of $\Rbar$ with no divisors lying in $R$ (the details and construction of these elements can be found in \cite{Kettinger2025}). Then:  
    $$[11(325+42\sqrt{-14})^{11}][11(325-42\sqrt{-14})^{11}]=11^2\cdot 19^{44}$$ is an element with elasticity at least $$\frac{2+44}{2}=23.$$ Since we know that $23\leq \rho(R)\leq\frac{D(\Cl(R))+1}{2}$, then necessarily $D(\Cl(R))=48$, so $\Cl(R)\cong\bZ_{48}$. Note that since $\ker(\tau)\cong \bZ_{12}$, $R$ must lie in case 2 of Theorem \ref{main result}, and thus $\rho(R)=24$.
\end{example}

    \bibliographystyle{plain}
    \bibliography{bibliography}

\begin{thebibliography}{10}

\bibitem{carlitz}
L.~Carlitz.
\newblock A characterization of algebraic number fields with class number two.
\newblock {\em Proceedings of the American Mathematical Society}, 11(3):391, 1960.

\bibitem{choi2024class}
Hyun~Seung Choi.
\newblock Class group and factorization in orders of a {PID}.
\newblock {\em Journal of Number Theory}, 265:226--269, 2024.

\bibitem{conrad}
Keith Conrad.
\newblock The conductor ideal of an order.
\newblock {\em Expository Paper}, 2019.
\newblock \url{https://kconrad.math.uconn.edu/blurbs/gradnumthy/conductor.pdf}.

\bibitem{radicalconductor}
James~Barker Coykendall and Grant Moles.
\newblock Elasticity in orders of an algebraic number field with radical conductor ideal and their rings of formal power series, 2025.

\bibitem{davenport1966midwestern}
H~Davenport.
\newblock Midwestern conference on group theory and number theory.
\newblock {\em Ohio State University}, 1966.

\bibitem{halter-koch}
Franz Halter-Koch.
\newblock Factorization of algebraic integers.
\newblock {\em Ber. Math. Stat. Sektion Forschung}, 191, 1983.

\bibitem{halter1995elasticity}
Franz Halter-Koch.
\newblock Elasticity of factorizations in atomic monoids and integral domains.
\newblock {\em Journal de th{\'e}orie des nombres de Bordeaux}, 7(2):367--385, 1995.

\bibitem{Kettinger2025}
Jared Kettinger.
\newblock A generalized davenport constant of the second kind, 2025.
\newblock \url{https://arxiv.org/abs/2506.13532}.

\bibitem{lmfdb}
The {LMFDB Collaboration}.
\newblock The {L}-functions and modular forms database.
\newblock \url{https://www.lmfdb.org}, 2024.
\newblock [Online; accessed 16 April 2024].

\bibitem{thesis}
Grant Moles.
\newblock The {HFD} property in orders of a number field.
\newblock {\em All Theses}, 3851, 2022.
\newblock \url{https://tigerprints.clemson.edu/all_theses/3851/}.

\bibitem{dissertation}
Grant Moles.
\newblock Relating elasticity and other multiplicative properties among orders in number fields and related rings.
\newblock {\em All Dissertations}, 3750, 2024.
\newblock \url{https://open.clemson.edu/all_dissertations/3750/}.

\bibitem{quadla}
Grant Moles.
\newblock {Associated Quadratic Orders}, March 2025.
\newblock \url{https://github.com/gramoles/Associated-Quadratic-Orders}.

\bibitem{narkiewicz}
Władysław Narkiewicz.
\newblock A note on elasticity of factorizations.
\newblock {\em Journal of Number Theory}, 51:46--47, 1995.

\bibitem{picavet}
Martine Picavet-L'Hermitte.
\newblock Some remarks on half-factorial orders.
\newblock {\em Rendiconti del Circolo Matematico di Palermo}, 52:297--307, 2003.

\bibitem{rago}
Balint Rago.
\newblock A characterization of half-factorial orders in algebraic number fields.
\newblock {\em arXiv preprint arXiv:2304.08099}, 2024.

\bibitem{valenza}
Robert~J. Valenza.
\newblock Elasticity of factorization in number fields.
\newblock {\em Journal of Number Theory}, 36:212--218, 1990.

\bibitem{zaks1}
Abraham Zaks.
\newblock Half factorial domains.
\newblock {\em Bulletin of the American Mathematical Society}, 82(5):721–723, Sep 1976.

\bibitem{zaks2}
Abraham Zaks.
\newblock Half-factorial-domains.
\newblock {\em Israel Journal of Mathematics}, 37:281--302, 1980.

\end{thebibliography}

\end{document}